\documentclass[10pt]{article}
\usepackage{amsmath, amsfonts, amsthm, amssymb, verbatim,dsfont,a4wide}
\usepackage{graphicx}
\usepackage[mathcal]{euscript}
\usepackage[applemac]{inputenc}
\usepackage{dsfont} 
\usepackage{mathrsfs}
\theoremstyle{plain}
        \newtheorem{theorem}{Theorem}[section]
        \newtheorem{proposition}[theorem]{Proposition}
        \newtheorem{lemma}[theorem]{Lemma}
        
\theoremstyle{definition}
        \newtheorem{definition}[theorem]{Definition}
        \newtheorem{remark}[theorem]{Remark}
\numberwithin{equation}{section}
\newcommand \be          {\begin{equation}}
\newcommand \ee        	{\end{equation}}
\newcommand \Dcal   	{\mathcal{D}} 
\newcommand \RR         	{\mathbb{R}} 
\newcommand \del        	\partial
\newcommand \eps       	\varepsilon
\newcommand \ubar     	{{\overline u}}
\newcommand \vbar       {{\overline v}}
\newcommand \Fbar       {{\overline F}}
\newcommand \wbar      {{\overline w}} 
\newcommand \loc        	{{\mathrm{loc}}}
\newcommand \ws     	{\mathrel{\mathop{\rightharpoonup}\limits^{*}}}
\newcommand \one 	{\mathds{1}}
\DeclareMathOperator    \sgn {sgn}
\DeclareMathOperator    \Real {Re}
\DeclareMathOperator    \Imag {Im}
\DeclareMathOperator \supp {supp} 
\begin{document}

\title{The Linear Stability of Shock Waves
 for the Nonlinear 
\\
Schr\"odinger--Inviscid Burgers System}

\author
{
Paulo Amorim$^1$, Jo\~ao-Paulo Dias$^1$,
\\
M\'ario Figueira$^1$, 
and
%
Philippe G. L{\large e}Floch$^2$
}  
\footnotetext[1]{Centro de Matem\'atica e Aplica\c c\~oes
Fundamentais, Universidade de Lisboa, Avenida Prof. Gama Pinto 2,
1649--003 Lisboa, Portugal. Email: {pamorim@ptmat.fc.ul.pt, dias@ptmat.fc.ul.pt, figueira@ptmat.fc.ul.pt}.}
\footnotetext[2]{
Laboratoire Jacques--Louis Lions 
\& Centre National de la Recherche Scientifique, Universit\'e Pierre et Marie Curie (Paris 6), 75252 Paris, France. 
Email : {contact@philippelefloch.org.} 
\newline  
\noindent
\textit{\ AMS Subject Class.} {35L65, 35L67, 35B35.} \hfill 
\hfill {\sl Published in J. Dyn. Diff. Equat. (2013). }
}  

\date{December 2012}

\maketitle

\begin{abstract} 
We investigate the coupling between the nonlinear Schr\"odinger equation and the inviscid Burgers equation, 
a system which models interactions between 
short and long waves, for instance in fluids. Well-posedness for the associated Cauchy problem remains a difficult
open problem, and we tackle it here via a linearization technique. Namely, we establish a linearized stability theorem 
for the Schr\"odinger--Burgers system, when the reference solution
is an entropy--satisfying shock wave to Burgers equation. Our proof is based on suitable energy estimates 
and on properties of hyperbolic equations with discontinuous coefficients. Numerical experiments support and expand our theoretical results. 
\end{abstract}


\section{Introduction} 
\label{sec-10}

\subsection{Background}

In \cite{Benney}, Benney introduced and studied several systems of partial differential equations 
which model the interaction of short waves and long waves arising in a variety of physical applications. These models are relevant, for instance, for the description of gravity waves in fluids or internal--surface waves. Several such systems 
were studied in the literature, 
especially systems that introduced a coupling between the nonlinear Schr\"odinger equation (for the short waves) and a general partial differential equation (for the long waves). In the literature, the long wave is often described by a linear transport equation or a nonlinear hyperbolic balance law
(which might also include dispersive or dissipative terms). More precisely, the general class of Benney's systems reads
\be
\label{50} 
\aligned 
i \, u_t + i c_1  u_x + u_{xx}
& = \alpha u\,v + \gamma |u|^2 u,
\\
v_t +  c_2 v_x + \mu v_{xxx} + \nu (v^2)_x 
& = \beta (|u|^2)_x,
\endaligned 
\ee
where $c_1,c_2, \alpha, \beta, \gamma, \mu$ and $\nu$ are real constants, 
and the unknown functions $u=u(t,x)$ and $v=v(t,x)$ represent the short waves and the long waves, respectively. For instance, 
Tsusumi and Hatano \cite{TsutsumiHatano1,TsutsumiHatano2} studied the coupling with a general linear equation in $v$, 
while Bekiranov, Ogawa, and Ponce  \cite{BOP1} studied the coupling with the Korteweg--de Vries (KdV) equation. 
More recently, there has been a renewed interest in the (particularly challenging) coupling with nonlinear balance laws, 
which motivates us here  
to consider the following \emph{nonlinear Schr\"odinger--inviscid Burgers} system (with $\eps>0$):   
\be
\label{100} 
\aligned
i u_t + u_{xx} & = vu - \eps |u|^2 u,
\\
v_t +  (v^2)_x & = \eps(|u^2|)_x. 
\endaligned  
\ee
This is clearly an important prototype for the class of systems \eqref{50}. 
In \cite{DFO1}, the existence of local--in--time and smooth solutions to this system was established by solving the initial value problem, 
but the possibility that smooth solutions blow--up for sufficiently large times was left open 
(although some partial results are available in \cite{ADFO}). 
On the other hand, the existence of {\sl global--in--time} solutions to the corresponding initial value problem 
is known \cite{DiasFigueira} only in the class of weak solutions, {\sl provided} the nonlinear flux $v^2$ in \eqref{100} is replaced by the cubic function $v^3$
(or, more generally, by any concave/convex function).

We recall also that Dias, Figueira, and Frid \cite{DFF} (cf.~also \cite{DiasFigueira,DiasFrid})
 treated a class of Benney--type systems in which an {\sl additional coupling} function was introduced in order to ``tame'' the nonlinear coupling between the two equations in \eqref{50}. These authors treated long waves modeled by, for instance, the Navier-Stokes equation or 
the $p-$system of nonlinear elasticity, and established a well-posedness theorem for the Cauchy problem.
Their theory encompasses only variants of the Schr\"odinger--Burgers system but, due to the presence of the additional coupling function in  \eqref{50}, 
the Schr\"odinger-Burgers system \eqref{100} of interest in the present work is not included.
Finally, we refer to \cite{AmorimDias} for the coupling with nonlinear viscoelasticity and to \cite{AmorimFigueira1,AmorimFigueira2,AmorimFigueira3,Chineses} for the numerical analysis of short waves--long waves problems. 


\subsection{Objective of this paper} 

Despite recent efforts on this subject, the existence of weak solutions to the Schr\"odinger--Burgers system remains an open problem.  In the present work, we shed some light on the systems \eqref{50} via a linearization approach and, specifically, we establish that solutions to  \eqref{100} are linearly stable, when the reference solution consists of an entropy--satisfying shock wave to the inviscid Burgers equation. Observe that the linearization of nonlinear hyperbolic equations (and systems) 
in a class of weak solutions was first investigated in LeFloch~\cite{CL,LeFloch0,LeFloch-book,LX}. More recently, this question was revisited by Godlewski and Raviart \cite{Raviart1} from the numerical standpoint and for systems in several space dimensions.  Concerning linear hyperbolic equations with discontinuous coefficients, arising below,  
we also refer the reader to~Bouchut and James~\cite{BouchutJames}, who introduced a notion of duality solutions for scalar hyperbolic equations. 

We adopt the following general strategy. Given  a system of partial differential equations and  a solution $U=U^\delta(t,x)$ which we assume to 
depend upon a (small)  parameter $\delta>0$, we seek for a formal expansion of that solution in the form
\be
\label{102}
\aligned 
U^\delta & = U^{(0)} + \delta \, U^{(1)} + \ldots,
\endaligned
\ee 
where $U^{(0)}$ is a chosen reference solution.
Formally at least, one can linearize the given (nonlinear) system around $U^{(0)}$ and derive 
a system for the first-order perturbation $U^{(1)}$. If we manage to establish an estimate of the form (on each compact time interval $[0,T]$, say)
\be
\label{103} 
\| U^{(1)}(t) \| \le G_T\big( \|U^{(1)}(0)\| \big) 
\ee
in a certain norm $\| \cdot \|$ and some continuous function $G_T:\RR_+ \to \RR_+$ vanishing at the origin, then we may expect that the nonlinear stability estimate 
$\| U^\delta(t) - U^{(0)}(t)\| \le G_T \big( \| U^\delta(0) - U^{(0)}(0)\| \big)$ would then follow.
This approach by linearization was adopted by LeFloch in~\cite{LeFloch-book}  in order to deal 
with solutions to nonlinear hyperbolic systems and estimate in the $L^1$ norm between two solutions was derived therein. This strategy was thus validated in this context and  led indeed to a proof of the $L^1$ continuous dependence of bounded variation solutions to nonlinear 
hyperbolic systems. 

The study of the full nonlinear problem \eqref{100} remains a very challenging open problem and, in the present work, 
we find it natural to focus on the linear stability issue which, as we will see, is quite involved. 
We will analyze solutions to a suitably linearized version of the Schr\"odinger--Burgers system \eqref{100}, 
when the reference solution consists of a shock wave (in $v$) and a smooth oscillatory traveling wave (in $u$). 
Our main result will be a rigorous stability statement of the form \eqref{103}.


\subsection{Outline of this paper}

To start off our analysis, we need to describe the solution that will serve as a reference state. For clarity in the presentation, we
consider first the limiting case $\eps=0$. It is not difficult to check that there exists a (stationary) traveling wave solution 
$$
\aligned
(u, v)(x) & = ( e^{ibt} r(x), \varphi(x)), \qquad x \in \RR, 
\\
\varphi(x) & = -\xi \sgn(x), \qquad \qquad x \in \RR, 
\endaligned
$$
in which $b$ is an arbitrary real parameter and $\varphi$ is a stationary shock wave  (normalized, from now on and without loss of generality, so that $\xi=1$) to Burgers equation,
 while the function $r \in W^{2,\infty}(\RR)$ is a solution to the ordinary differential equation  
\be
\label{170}
\aligned
- r(x)'' = -b r(x) + r(x) \sgn x. 
\endaligned
\ee
We look for solutions that are continuous with continuous derivative at $x=0$, so that, with obvious notation,
\[
\aligned
r(x) = \one_{x>0} r_+(x) + \one_{x<0} r_-(x),
\endaligned
\]
and 
$r_+(0) = r_-(0)$ and $r'_+(0) = r'_-(0)$. Elementary methods lead to us, if $b\le -1$, to 
\be
\label{175}
\left\{\aligned
& r_+(x) = C \, \sqrt{\frac{|1+b|}{1-b}} \sin\big( \sqrt{1-b} \, x \big) + A \, \cos\big( \sqrt{1-b} \, x \big),
\\
& r_-(x) = C \, \sin\big( \sqrt{|1+b|} \, x \big) + A \, \cos\big( \sqrt{|1+b|} \, x\big),
\\
& r(0) = A, \qquad r'(0) = C \sqrt{|1+b|},
\endaligned
\right.
\ee
 and, if $-1 < b< 1$, to  
\be
\label{180}
\left\{\aligned
& r_+(x) = A \, \sqrt{\frac{b+1}{1-b}} \sin\big( \sqrt{1-b} \, x \big) + A \, \cos\big( \sqrt{1-b} \, x \big),
\\
& r_-(x) = A e^{\sqrt{b+1}x},
\\
& r(0) = A, \qquad r'(0) = A \, \sqrt{1+b},
\endaligned
\right.
\ee
where $A, C$ are arbitrary constants.

More generally, when $\eps>0$, the corresponding solution $r_\eps(x)$ is given by
\be
\label{150}
\aligned
&r_\eps''(x) = b \, r_\eps(x) - r_\eps(x)\sgn x \, \sqrt{\eps r^2_\eps(x) + 1} - \eps \, r_\eps^3(x),
\\
& r_\eps(0) = r(0), \qquad r'_\eps(0) = r'(0),
\endaligned
\ee
and 
\be
\label{160}
\aligned
\varphi_\eps (x) = -\sgn x \, \sqrt{\eps r^2_\eps(x) + 1} \in L^\infty(\RR). 
\endaligned
\ee
No closed formula is now available but standard arguments for ordinary differential equations 
(i.e.~Picard's theorem and an energy estimate) 
yield the existence and uniqueness of a global solution $r_\eps \in W^{2,\infty}(\RR)$ of \eqref{150}. 
An argument of continuity with respect to parameters shows that, as $\eps \to 0$, 
one has $r_\eps \to r$ in $C^1_\loc(\RR)$ (uniform convergence of first--order derivatives), 
where $r$ is defined by \eqref{175}--\eqref{180}. 

Obviously, the choice of a minus sign in \eqref{160} is not arbitrary. Indeed, one could attempt to 
perform the foregoing analysis with the (monotone increasing) weak solution 
$\varphi_\eps (x) = \sgn x \sqrt{\eps r^2_\eps(x) + 1}$. However, as we shall see, the choice \eqref{160} 
corresponds to an entropy--satisfying shock and guarantees uniqueness to the linearized problem. (Infinitely many solutions would otherwise be available.) Hence, the system \eqref{100} may be viewed as a diffusive regularization of Burgers equation which 
properly takes into account the small-scale effects, and thus allows us to consider the whole model as a relevant physical model.

An outline of the paper follows. In Section~\ref{sec-200}, we consider the weakly coupled problem obtained by setting the parameter $\eps$ to zero.
This regime is physically interesting, since (as was argued in \cite{DFF}) it reflects the fact that a shock wave is a macroscopic phenomenon 
while the Schr\"odinger equation models a microscopic one. 
One main difficulty stems from the fact that the linearized system involves a transport equation with discontinuous
coefficient. In the case $\eps =0$ the linearized hyperbolic equation can be solved
explicitly and yields a singular solution containing a Dirac measure \cite{LeFloch0}.
Our main task is to analyze the interaction between this measure--solution and the solution to the Schr\"odinger equation.

In Section~\ref{sec-400}, we consider the full problem with $\eps >0$ and we establish our main stability result. The analysis 
of the linearized system in this general regime is 
made challenging by the fact that the Schr\"odinger equation now has a nonlinear right-hand side coupled to the hyperbolic equation, and
it is convenient to rely on a change of variables introduced earlier in \cite{DFO1,AmorimDias} in order 
to derive a formulation amenable to Kato's semi--group method \cite{Kato2}. Moreover, the linearized hyperbolic equation
still involves the product of a measure by a discontinuous function and can no longer be solved explicitly. It is
more delicate to define a suitable notion of solution. Finally, in Section~\ref{sec-600}, we present some numerical 
experiments and illustrate our theoretical results.


\section{Linearized stability for a weakly coupled system}
\label{sec-200}

\subsection{Stability statement}

Our first objective is the linearized analysis of the system \eqref{100} when the coupling parameter $\eps$ is set to zero. Thus, we consider the weakly coupled system ($t>0, x \in \RR$) 
\be
\label{500} 
\aligned
&i u_t + u_{xx} = vu, 
\\
&v_t +  (v^2)_x = 0,
\endaligned 
\ee
We linearize this system around the solution $\big( e^{ibt} r, \varphi \big)$ described in the introduction, especially in \eqref{175} and \eqref{180}. We make the substitution
\[
v \mapsto \varphi + \delta v,
\qquad
\quad
u \mapsto e^{ibt}r + \delta e^{ibt} u 
\]
for some small parameter $\delta>0$. By identifying first--order terms in $\delta$, we find 
\be
\label{550} 
\aligned
&i u_t + u_{xx} = (\varphi + b)u + vr, 
\\
&v_t +  2(\varphi v)_x = 0, 
\endaligned  
\ee
which we supplement with the initial data
\be
\label{555}
\aligned
u(x,0) = u_0 (x) \in H^1(\RR), \qquad v(x,0) = v_0(x) \in H^1(\RR).
\endaligned
\ee
(The notion of weak solution used in the following theorem will be clarified shortly.)  

\begin{theorem}[Linearized stability for a weakly coupled system]
\label{thm-100}
Let $r=r(x)$ be defined by \eqref{175} or \eqref{180} and $\varphi = -\sgn(x).$ Then, the solution 
$(e^{ibt} r, \varphi)$ to system \eqref{500} is linearly stable in the following sense:
for any $T>0$ and any initial data \eqref{555}, the system \eqref{550} admits a unique weak solution $(u,v)= (u(t,x), v(t,x))$ 
satisfying   
\be
\label{560}
\aligned
\| (u , v) \|_{ L^\infty(0,T; H^1 \times H^{-1})} \le G_T(\| (u(0), v(0)) \|_{H^1 \times (L^1 \cap L^2)}),
\endaligned
\ee
where $G_T: \RR_+ \to \RR_+$ is a continuous function vanishing at the origin.  
\end{theorem}

The solution $v$ in Theorem~\ref{thm-100} was first introduced in \cite{LeFloch0} and can also be interpreted 
as the unique duality solution to the second equation in \eqref{550} in the sense of~\cite{BouchutJames}.  


\subsection{Proof of stability} 
\label{sec-300}

We now discuss the relevant notion of weak solutions and we establish the linearized stability property of interest. 
We need to establish existence and uniqueness for the problem \eqref{550}--\eqref{555}
and derive the bound \eqref{560}. 
First of all, we observe that the second  equation in \eqref{550} admits an explicit solution 
(determined in \cite{LeFloch0}, even for general Riemann data):
\be
\label{600}
\aligned
v(t,x) = \widetilde v(t,x) + \Psi(t) \delta_{\Sigma_T},
\endaligned
\ee
where 
$$
\widetilde v(t,x) = \left\{
\aligned
&v_0(x-2t), && x<0,
\\
&v_0(x+2t), && x>0,
\endaligned
\right.
\qquad
\quad
\Psi(t) =  \int_{-2t}^{2t} v_0(x) \,dx,
$$ 
and $\delta_{\Sigma_T}$ denotes the Dirac measure concentrated on the line segment $\Sigma_T = \{ (0,t) : t \in(0,T) \}$.
Thus, if $\phi$ is a test function in $\RR^2$, we set  
$\langle \delta_{\Sigma_T}, \phi \rangle = \int_0^T \phi(0,t) \,dt$. 
In all that follows, and since the final time $T>0$ is arbitrary but fixed, we simply use $\Sigma$ to denote $\Sigma_T$.

Some remarks about \eqref{600} are in order. First, observe that $v$ is clearly a solution to \eqref{500} outside the line $\Sigma_T$.
Second, observe that replacing the expression of $v$ in \eqref{500} leads to a product
 of the discontinuous function $2\varphi$ by the distribution $\delta_\Sigma$, which has no meaning in distribution theory. 
Therefore, in agreement with the definition introduced by LeFloch~\cite{LeFloch0}, one defines 
$2\varphi \delta_\Sigma = \sigma$, where $\sigma$ is the speed of propagation of the shock $\varphi$ which coincides with the direction of $\Sigma$ (in this case, $\sigma=0$).
This gives 
\[
\aligned
v_t + 2(\varphi v)_x = \widetilde v+ 2(\varphi \widetilde v)_x + \Psi'(t) \delta_\Sigma + \Psi(t) \big( \del_t + \sigma \del_x\big) \delta_\Sigma,
\endaligned
\]
so that the definition of the product $2\varphi \delta_\Sigma$ merely translates the natural fact that $\delta_\Sigma$ is invariant in the direction of $\Sigma$. With this in mind, it is easy to see that $v_t + 2(\varphi v)_x = 0$.

Let us now show the estimate 
\be
\label{620}
\aligned
\| v \|_{ L^\infty(\RR_+; H^{-1}(\RR))} \le c (\| v_0\|_{ L^1(\RR)} + \| v_0\|_{ L^2(\RR)}),
\endaligned
\ee
which will establish the linearized stability for the second equation in \eqref{500}. First, it is obvious that
$
\| \widetilde v \| \le \| v_0\|_{ L^2}.
$
Next, given $\phi=\phi(t,x) \in L^1(0,T; H^1)$, we can write 
\[
\aligned
| \langle \Psi(t) \delta_\Sigma , \phi \rangle | &= | \langle \delta_\Sigma , \Psi(t) \phi(t,x) \rangle|
\\
& = \Big| \int_0^T \int_{-2t}^{2t} v_0(x) \phi(0,t) \,dx \,dt \Big|
\leq \| v_0 \|_{L^1} \int_0^T \| \phi(t) \|_{L^\infty(\RR) } \,dt.
\endaligned
\]
Since $H^1(\RR) $ is continously embedded in $ L^\infty(\RR)$, we obtain
$$
| \langle \Psi(t) \delta_\Sigma, \phi \rangle | \le c\| v_0\|_{L^1} \| \phi\|_{L^1(\RR_+; H^1)},
$$
whence \eqref{620}.

We now turn our attention to the Schr\"odinger equation in \eqref{550}. In order to establish the existence of a solution, we regularize the singular term $v \, r$ and, for each $n>0$,
we consider the problem
\be
\label{700}
\aligned
&i u^n_t + u^n_{xx} = V u^n + h^n r,
\\
&u^n(x,0) = u^n_0(x),
\endaligned
\ee
where 
\be
\label{710}
\aligned
h^n(t,x) := \widetilde v(t,x) + \Psi(t) n\rho(nx),
\endaligned
\ee
where $\rho$ is a standard mollifier, $V(x) = \varphi + b$, $\varphi = -\sgn x$, and $\widetilde v$ is given in \eqref{600}.
In addition, $u^n_0 \in H^2(\RR)$ is some regularization of the data $u_0$, but 
our final estimate will actually hold with $u_0 \in H^1$, only.

Standard techniques (cf.~for instance \cite{Cazenave}) about the semi--group generated by $i( \Delta - V)$ with a potential $V\in L^\infty$ imply that, for each $n$, the problem \eqref{700} with $u^n_0 \in H^2$ admits a unique solution 
$u^n \in C^1 \big( [0,T]; L^2 \big) \cap C \big( [0,T]; H^2 \big)$
for all $T>0$, provided that (for fixed $n>0$) $h^n \in W^{1,1}(0,T; L^2(\RR))$.
This latter property follows easily from the fact that
$\del_t h_n = \del_t \widetilde v + 2(v_0 (2t) + v_0(-2t)) n \rho(nx)$, 
so that (for fixed $n$)
\[
\aligned
\| \del_t h_n\|_{ L^1(0,T; L^2)} \le c \| v_0 \|_{L^1(0,T; H^1)} .
\endaligned
\]

In oder to justify the passage to the limit $n\to\infty$ in the equation \eqref{700}, we need first obtain a (uniform in $n$) estimate for $h^n$, namely
\be
\label{720}
\| h_n\|_{ L^\infty(0,T; H^{-1}(\RR))} \le c (\| v_0\|_{L^1(\RR)} + \| v_0\|_{ L^2(\RR)})
\ee
for some uniform $c$. Given $\phi=\phi(t,x) \in L^1(0,T; H^1)$, by neglecting the term involving $\widetilde v$
(which was already dealt with) and by using the unit integral property of the mollifier, we find 
\[
\aligned
| \langle \Psi(t) n\rho(nx), \phi(t,x) \rangle_{H^{-1}\times H^1}| &\le \int_0^T \int_{\RR} \int_{-2t}^{2t} |v_0 (s)| \,ds\, n\rho(nx) |\phi(t,x)| \,dx \,dt
\\
&\le \|v_0 \|_{L^1} \int_0^T \| \phi\|_{L^\infty(\RR)} \,dt
\leq c \|v_0 \|_{L^1} \int_0^T \| \phi\|_{H^1} \,dt,
\endaligned
\]
from which \eqref{720} follows.

Next, we derive an estimate of the $L^2$-norm of $u^n$. Consider the equation \eqref{700},
 multiply it by $\ubar^n$, and take the imaginary part: 
$\Imag (i \ubar u_t + \ubar u_xx ) = \Imag (\ubar r h_n)$, 
where, from now on, we omit the superscript $n$. Since $\Imag i\ubar u_t = \frac12 \del_t |u|^2$, by integrating 
over $\RR$, we find 
\[
\frac12 \frac d{dt} \int_{\RR} |u|^2 \,dx  = \int_{\RR} r \ubar h_n \,dx,
\]
where we have used $\Imag \int_{\RR} |u_x|^2 =0$. Integrating over $(0,T)$ yields
\[
\aligned
\frac12 \int_{\RR} |u(t)|^2 \,dx &\le c + \int_0^t \Big| \int_{\RR} r\ubar h_n \,dx \Big| \,d\tau
 \le c + c \int_0^t |\langle h_n, \ubar \rangle_{H^{-1}\times H^1}| \,d\tau
\endaligned
\]
which, by \eqref{720}, gives
\be
\label{750}
\aligned
\| u^n\|^2_{L^2} \le c +c \int_0^t \| u^n\|_{H^1} \, d\tau.
\endaligned
\ee

We now estimate the gradient of $u^n$, as follows. Multiply the equation \eqref{700} by $\ubar_t$, take the real part, and integrate over $\RR$ in order to obtain 
\[
\aligned
\frac12 \frac d{dt} \int_{\RR} | u_x|^2 + V |u^2| \,dx & = -\int_{\RR} r h_n \Real (u) \,dx
\\
& = -\frac d{dt} \int_{\RR} r h_n \Real (u) \,dx + \int_{\RR} r \del_t h_n \Real (u) \,dx.
\endaligned
\] 
Then, an integration over $(0,T)$ gives (since $r$ is bounded)
\be
\label{760}
\aligned
& \frac12 \int_{\RR} | u_x|^2 \,dx + \frac12 \int_{\RR} V |u^2| \,dx  
\\
&\le c + c \Big| \int_{\RR}  h_n \Real (u) \,dx \Big|
+ c\int_0^t \Big| \int_{\RR} \del_t h_n \Real (u) \,dx \Big| \,d\tau.
\endaligned
\ee
We now estimate the terms on the right-hand side. First, using \eqref{720} and a weighted Young inequality, one has 
\be
\label{770}
\aligned
c \, \Big| \int_{\RR}  h_n \Real (u) \,dx \Big| &\le c (\| v_0\|_{L^1}^2 + \| v_0\|_{ L^2}^2) + \frac14 \|u\|_{H^1}
\endaligned
\ee
and, in view of \eqref{710}, $\del_t h_n = \del_t \widetilde v + 2(v_0 (2t) + v_0(-2t)) n \rho(nx)$, we find
\[
\aligned
\int_0^t \Big| \int_{\RR} \del_t h_n \Real (u) \,dx \Big| \,d\tau 
& \le \int_0^t \int_{\RR} |\widetilde v_t u| \,dx
+ c \int_{\RR} |v_0 (2t) + v_0(-2t)| n \rho(nx) |u| \,dx \,d\tau
\\
& \le c \int_0^t \|v_0\|_{H^1} \|u\|_{L^2} + \|v_0 u\|_{L^\infty} \,d\tau.
\endaligned
\]
Since $H^1 \subset L^\infty$, it follows that 
\be
\label{780}
\aligned
&\int_0^t \Big| \int_{\RR} \del_t h_n \Real (u) \,dx \Big| \,d\tau &\le \int_0^t\|v_0\|_{H^1}\|u\|_{H^1} \,d\tau.
\endaligned
\ee
Next, using \eqref{770} and \eqref{780} in \eqref{760} gives
\[
\aligned
\frac14 \| u^n_x \|^2_{L^2} \le c + c \| u^n\|^2_{L^2} + c \int_0^t \|u^n\|_{H^1} \,d\tau,
\endaligned
\]
which, along with \eqref{750}, implies 
$\| u^n\|_{H^1}^2 \le c +c \int_0^t \| u^n\|_{H^1} \,d\tau$. 
An application of Gronwall's lemma then leads us to 
$\| u^n\|_{H^1} \le c(t)$ 
for some continuous function of $t$ that is independent of $n$.

Thus, there exists $u \in L^\infty(0,T; H^1)$ such that as $n\to \infty$
(after the extraction of a subsequence) 
\[
\aligned
&u^n \ws u \quad L^\infty(0,T; H^1) \quad \text{ weak }*,
\\
&u^n_{xx} \ws u_{xx} \quad L^\infty(0,T; H^{-1}) \quad  \text{ weak }*,
\\
&Vu^n \ws Vu \quad L^\infty(0,T; L^2) \quad \text{ weak }*. 
\endaligned
\]
Furthermore, since $h^n \to v$ in the sense of distributions (with $v$ given by \eqref{600}), 
it follows from \eqref{720} that 
$h^n \, r \ws v \, r$ in $L^\infty(0,T; H^{-1})$ weak $*$.  
The equation for $u^n$ then yields 
$u^n_t \ws u_t$ in $L^\infty(0,T; H^{-1})$ weak $*$  
and we conclude that $u$ satisfies the first equation in \eqref{550}. Uniqueness is a consequence of the linearity of the equation.
It is possible to deduce additional regularity. Indeed, since 
$u \in L^2(0,T; H^1)$ and $u_t \in L^2(0,T; H^{-1})$, 
then, (cf.~\cite[Chap. 2]{Lions}) we can assume
$u\in C([0,T];L^2)$, which completes the proof of Theorem~\ref{thm-100}.


\section{Linearized stability for the Schr\"odinger--Burgers system}
\label{sec-400}

\subsection{Stability statement}

We now turn to the linearization of the full system \eqref{100} (with $\eps>0$) and consider the reference solution
 $(e^{ibt} r_\eps, \varphi_\eps)$ given in Section~\ref{sec-10}; cf.~\eqref{150},\eqref{160}.
After a routine calculation, we arrive at the following
\emph{linearized Schr\"odinger--Burgers system}
\be
\label{200} 
\aligned
&i u_t + u_{xx} = (\varphi_\eps + b -2 \eps \, r_\eps^2)u - \eps r_\eps^2 \, \ubar + v \, r_\eps, 
\\
&v_t +  2(\varphi_\eps v)_x = 2 \eps \, \Real( r_\eps u)_x,
\endaligned  
\ee
supplementd with the initial data
\be
\label{201}
\aligned
u(x,0) = u_0 (x) \in H^1(\RR), \qquad v(x,0) = v_0(x) \in H^1(\RR).
\endaligned
\ee

We introduce the following notion of weak solution.   

\begin{definition}
\label{def-200}
Consider the reference functions $\varphi_\eps$ and $r_\eps$ in \eqref{150}--\eqref{160} for some $\eps>0$.  
A pair $(u,v) \in L^\infty(0,T; H^1 \times H^{-1}(\RR))$ is called 
a {\rm weak solution} to the Cauchy problem \eqref{200}-\eqref{201}
if the following conditions hold: 
\begin{enumerate}
\item $u$ is a weak solution to the Schr\"odinger equation
\[
\aligned
&i u_t + u_{xx} = (\varphi_\eps + b -2 \eps r_\eps^2)u - \eps r_\eps^2 \ubar + vr_\eps. 
\endaligned
\]
\item $v \in L^\infty(0,T; H^{-1}(\RR)) $ has the form
\be
\label{220}
\aligned
v(t,x) = \widetilde v(t,x) + \Psi(t) \delta_\Sigma,
\endaligned
\ee
where $\widetilde v \in L^\infty(0,T; L^2(\RR))$, 
$\Sigma := \{(t,x) : x=0, t\in (0,T)\}$,
$\delta_\Sigma$ is the Dirac measure on $\Sigma$, 
$\Psi(t) \in L^\infty(0,T)$, and in the sense of distributions on $\RR \times (0,T)$
\be
\label{230}
\aligned
v_t +  2(\varphi_\eps v)_x = 2 \eps \Real( r_\eps u)_x,
\endaligned
\ee
with $\varphi_\eps \delta_\Sigma := 0$. 

\item
Finally, the initial data are assumed in a weak sense, that is, setting $u_0$ and $v_0$ in the integral equations involving appropriate test functions,
in the usual definition of weak solution to the Cauchy problem.
\end{enumerate}
\end{definition}

This definition agrees with the one in Section~\ref{sec-300} and  
(as far as the hyperbolic equation is concerned) generalizes the one in~\cite{LeFloch0} for hyperbolic equations (but with a non--vanishing right--hand side). 

\begin{remark} Of course, under the assumption that $v$ has the form \eqref{220} one can directly 
rewrite the equation for $v$ as 
$v_t +  2(\varphi_\eps \widetilde v)_x = 2 \eps \Real( r_\eps u)_x$, so that the condition 
 $\varphi_\eps \delta_\Sigma := 0$ need not be stated explicitly. Furthermore, 
from our proof of Theorem~\ref{thm-200} below it follows that a solution $v$ satisfying Definition~\ref{def-200} 
is actually the unique duality solution (in the sense of \cite{BouchutJames}) to the second equation of \eqref{200}
(as follows easily from Theorem~4.3.2 therein, by taking the source--term into account). 
\end{remark}

We are now in a position to state the main result of this paper. 

\begin{theorem}[Linearized stability property for the nonlinear Schr\"odinger--inviscid Burgers system] 
\label{thm-200}
Consider the reference functions $r_\eps=r_\eps(x)$ and $\varphi_\eps=\varphi_\eps(x)$ given in Section~\ref{sec-10}; cf.~\eqref{175}--\eqref{180}. Then, the solution $(e^{ibt} r_\eps, \varphi_\eps)$ to system \eqref{100} is linearly stable in the sense that, for any $T>0$, and any initial data \eqref{555}, the system \eqref{200} admits 
a unique weak solution $(u,v) = (u(t,x), v(t,x))$ in the sense of 
Definition~\ref{def-200} with 
\be
\label{250}
\| (u , v) \|_{ L^\infty(0,T; H^1 \times H^{-1})} \le G_T\big( \| (u(0), v(0)) \|_{H^1 \times H^1}\big), 
\ee
where $G_T:\RR_+ \to \RR_+$ is a continuous function  vanishing at the origin. 
\end{theorem}


\subsection{Proof of stability} 
\label{sec-500}

This section contains a proof of Theorem~\ref{thm-200}, that is, a proof of well posedness for the problem \eqref{200}--\eqref{201}. To this end, it is necessary to carefully study the coupling between the two equations
and deal with the lack of regularity of the coefficient of the second equation. We can first establish a well-posedness result for a regularized problem with smooth coefficients, as now stated. 

\begin{proposition}
\label{prop-400}
Let $a_j=a_j(x)$ ($j=1,\dots,5$) be real functions in $W^{2,\infty}(\RR)$ such that $a_5 =  a_3$ up to a multiplicative constant, and fix any data $u_0 \in H^2(\RR)$ and $v_0 \in H^1(\RR)$. Then, the problem
\be
\label{850} 
\aligned
&i u_t + u_{xx} = a_1 u + a_2 \ubar + a_3 v,
\\
& v_t + (a_4 v)_x = (a_5 \Real (u))_x,
\endaligned 
\ee
with prescribed initial data $(u(0),v(0)) = (u_0, v_0)$, admits a unique solution satisfying
\[
(u, v) \in C([0,+\infty[; H^2) \cap C^1( [0,+\infty[; L^2) \times C([0,+\infty[; H^1) \cap C^1([0,+\infty[; L^2). 
\]
\end{proposition}

\proof In order to construct the local--in--time (strong) solution to \eqref{850}, we follow the technique in \cite{Oliveira,DFO1}
and introduce an auxiliary system with non--local source, which can be tackled via  Kato's theory \cite{Kato2}. 
This is necessary in order to write the system \eqref{850} without derivative loss (see \cite{DFO1} for details). In other words, 
by refering to \cite{DFO1} for the motivation, we consider the linear system
\be
\label{900} 
\aligned
&i F_t + F_{xx} = a_1 F + a_2 \Fbar - a_3 a_4 v_x - a_3 a'_4 v - a_3 (a_5 \Real \widetilde u)_x, 
\\
& v_t + a_4 v_x + a_4' v = (a_5 \Real \widetilde u)_x, 
\endaligned 
\ee
where $\overline{F}$ is the complex conjugate of $F$ and
\be
\label{905}
\aligned
&u(t,x) = u_0(x) + \int_0^t F(x,s) \,ds,
\\
&\widetilde u(t,x) = (\Delta - 1)^{-1} \big( (a_1 -1) u + a_2 \ubar + a_3 v - i F\big)
\endaligned
\ee
with initial data 
\be
\label{910}
F(\cdot,0) = F_0 \in L^2(\RR), \qquad v(\cdot, 0) = v_0 \in H^1(\RR).
\ee
Once we have a local solution
\be
\label{920}
\aligned
F \in C([0,T]; L^2) \cap C^1([0,T]; H^{-2}), \quad v \in C([0,T]; H^1) \cap C([0,T]; L^2)
\endaligned
\ee
to the problem \eqref{900}, for some small $T>0,$ one can argue (as in \cite[Lemma 2.1]{DFO1}) and show that $(u,v)$ given by \eqref{900}--\eqref{905}
is actually the desired local solution. 

At this stage, we only sketch the argument since it is quite similar to the one in \cite{AmorimDias,DFO1}. First, we write \eqref{900} as a system of three equations, by decomposing $F$ into its real and 
imaginary parts, which allows us to obtain the abstract form
\be
\label{930} 
\aligned
& U_t + A(U) U = g(t,U), \qquad \quad 
U(\cdot, 0) = U_0,
\endaligned 
\ee
with $U = (\Real F, \Imag F, v)$ and some initial data $U_0$.
The key point is to decompose the operator 
\[
A(U) = 
\left[ 
\begin{array}{ccc} 
0 & \Delta & 0  
\\
-\Delta& 0 & 0 
\\ 
0 & 0 & a_4 \del_x + a'_4
\end{array} 
\right]
\]
in the form $S A(U) S^{-1} = A(U) + B(U)$ for some operator $B$ (again, see~\cite{AmorimDias,DFO1} for a similar
statement). In the present setting, we obtain such a decomposition by setting 
\[
S := 
\left[ 
\begin{array}{ccc} 
1-\Delta & 0& 0  
\\
0 & 1-\Delta & 0 
\\ 
0 & 0 & (1-\Delta)^{1/2}
\end{array} 
\right].
\]
Note that $S:Y \to X$ is an isomorphism, provided $Y := (L^2)^2 \times H^1$ and $X := (H^{-2})^2 \times L^2$.
The relevant properties satisfied by $S$ (in particular about $(1-\Delta)^{1/2}$) can be found in \cite[Section 8]{Kato1}.
Observe that the right--hand side of \eqref{930} is linear in $U$, so that it is straighforward 
to derive the necessary estimates
for the source $g$ and we may finally apply \cite[Theorem 6]{Kato2} and conclude with the existence of a unique
$F,v$ satisfying \eqref{920}. 

We have thus established a local existence result for the Cauchy problem \eqref{850}. In order to show that the solution exists for all time, we must derive certain \emph{a priori} estimates. We begin with an energy inequality satisfied by solutions 
to \eqref{850}. We multiply the first equation in \eqref{850} by $\ubar_t$, take the real
part, and integrate over $\RR$, and obtain 
\be
\label{940}
\aligned
\frac d{dt} \int_{\RR} | u_x|^2 \,dx + \frac d{dt} \int_{\RR} a_1 | u|^2 \,dx  + \frac d{dt} \Real \int_{\RR} a_2 u^2 \,dx = - 2 \Real \int_{\RR} a_3 v u_t \,dx.
\endaligned
\ee
By assumption $a_5 = c a_3$ for some constant $c$, so we can use the second equation in \eqref{850} and find  
\[
\aligned
\Real \int_{\RR} a_3 v u_t \,dx &= \Real \frac d{dt} \int_{\RR}a_3 vu \,dx - \Real \int_{\RR} a_3 v_t u \,dx
\\
&=  \Real \frac d{dt} \int_{\RR}a_3 vu \,dx + \Real \int_{\RR} a_3 u (a_4 v)_x - c\int_{\RR}a_3 \Real (u) (a_3 \Real (u))_x \,dx
\\
&=  \Real \frac d{dt} \int_{\RR}a_3 vu \,dx - \Real \int_{\RR} (a_3 u)_x a_4 v \,dx
\\
&=  \Real \frac d{dt} \int_{\RR}a_3 vu \,dx - \Real \frac d{dt} \int_{\RR} a_4 v^2 \,dx.
\endaligned
\]
Combining this result with \eqref{940}, we arrive at the energy conservation property
\be
\label{950}
\frac d{dt} \int_{\RR} \Big( |u_x|^2 + a_1 |u|^2 + \Real a_2 u^2 + 2 \Real a_3 vu - a_4 v^2 \Big) \, dx = 0
\ee
and an immediate consequence is thus 
\be
\label{960}
\aligned
\| u_x \|_2^2 \le c + c\| v\|_2^2 + c\| u\|_2^2.
\endaligned
\ee

Now, we multiply the second equation in \eqref{850} by $v$,
 integrate over $\RR$, and finally apply Gronwall's lemma (assuming that $v$ is  
sufficiently smooth, without loss of generality as can be checked by an elementary smoothing argument): 
$
\| v\|_2^2 \le c + c(t) \int_0^t \| u_x\|_2^2 + \|u \|_2^2 \,ds. 
$
Along with the estimate \eqref{960} and again by applying Gronwall's lemma, this gives 
\be
\label{965}
\aligned
\| v\|_2^2  &\le c(t) + c(t) \int_0^t \|u\|_2^2 \,ds. 
\endaligned
\ee
Next, multiplying the first equation in \eqref{850} by $\ubar$, taking the imaginary part, and integrating in $\RR$, we obtain 
\be
\label{966}
\aligned
\frac 12 \frac d{dt} \int_{\RR} |u|^2 \,dx = - \Imag \int_{\RR} a_2 (\ubar)^2 \,dx+ \Imag \int_{\RR} a_3 v \ubar \,dx.
\endaligned
\ee
Thus, an integration over $(0,t)$, Gronwall's lemma, and \eqref{965} yield
\[
\| u \|_2^2  \le c(t) + c(t) \int_0^t \|v^2\|_2^2 \,ds 
  \le c(t) + c(t) \int_0^t \|u\|_2^2 \,ds 
\]
and, with Gronwall's lemma again, 
$
\| u \|_2^2 \le c(t).
$
From \eqref{965}, we see that $\| v\|_2^2 \le c(t)$ and thus, in combination with \eqref{960}, we find
\be
\label{970}
\| u_x\|_2^2 + \| u\|_2^2 + \|v\|_2^2 \le c(t)
\ee
for some continuous function $c=c(t)$.

Next, we need to estimate $v_x, v_t, u_t $ and $u_{xx}$. First, differentiate the second equation in \eqref{850} with respect to $x$, 
multiply by $v_x$ and integrate over $\RR$. We obtain after some trivial integration by parts
\[
\aligned
\frac12 \frac d{dt}\int_{\RR} (v_x)^2 \,dx + \int_{\RR} a_4'' v v_x \,dx + \frac32 \int_{\RR} a_4'(v_x)^2 \,dx = \int_{\RR} (a_5 \Real (u))_{xx} v_x.
\endaligned
\]
Using Young's inequality and integrating in $(0,t)$ gives
\[
\aligned
\| v_x (t)\|_2^2 \le \| v_x(0)\|_2^2 + c \int_0^t \|u_{xx}\|_2^2 + \|v_x\|_2^2 \,dt
\endaligned
\]
and so, by Gronwall's lemma,
\be
\label{972}
\aligned
\| v_x (t)\|_2^2 \le c( 1 + t  \|u_{xx}\|_2^2) e^{ct}.
\endaligned
\ee
Next, the first equation in \eqref{850} and \eqref{970} give immediately
\be
\label{973}
\aligned
\| u_{xx} \|_2^2 \le \| u_t \|_2^2 + c(t).
\endaligned
\ee
Now, we take the time derivative of the first equation, multiply by $\ubar_t$, take the imaginary part and integrate in $\RR$ to get 
\be
\label{973a}
\aligned
\frac12\frac d{dt} \|u_t\|_2^2 \le c \|u_t\|_2^2 + c \|v_t\|_2^2. 
\endaligned
\ee
From the first equation, \eqref{970} and \eqref{972} we get
\[
\aligned
\|v_t\|_2^2 \le c \| v_x\|_2^2 + c(t) \le c(t)(1 + \|u_{xx}\|_2^2 ),
\endaligned
\]
and so \eqref{973} and \eqref{973a} give
\be
\label{974}
\aligned
\frac d{dt} \|u_t\|_2^2 \le c(t)( 1+ \|u_t\|_2^2).
\endaligned
\ee
Applying Gronwall's lemma, we conclude that $\|u_t\|_2^2 \le c(t)$. This estimate and \eqref{972}--\eqref{974} together give
\be
\label{975}
\aligned
\| v_x\|_2^2 + \| v_t\|_2^2 +  \|u_t\|_2^2 + \| u_{xx}\|_2^2 \le c(t)
\endaligned
\ee
for some continuous function $c(t)$. In view of these uniform estimates, it is now standard to verify that the local-in-time solution
is actually defined for all $t>0$. This concludes the proof of Proposition~\ref{prop-400}.
\endproof

To be able to apply Proposition \ref{prop-400}, we now introduce a regularized version of the linearized system \eqref{200}.  
For $\delta >0$, we define $\varphi^\delta_\eps$ from $\varphi_\eps$ by convolution of the sign function with a standard mollifier. (In particular,  $\varphi^\delta_\eps(0) = 0$.) 
Using Proposition~\ref{prop-400}, we obtain a solution 
$(u^\delta, v^\delta)$ in $C([0,T]; H^2) \cap C^1( [0,T]; L^2) \times C([0,T]; H^1) \cap C^1([0,T]; L^2)$ to 
the Cauchy problem ($\eps$ being omited)
\be
\label{980} 
\aligned
&i \, u_t + u_{xx} = (\varphi^\delta + b -2 \eps r^2) \, u - \eps r^2 \ubar + v \, r, 
\\
&v_t +  2(\varphi^\delta \, v)_x = 2 \eps \, \Real( r u)_x,
\endaligned 
\ee
with initial data
\be
\label{990}
\aligned
u^\delta(x,0) = u^\delta_0 (x) \in H^2(\RR), \qquad v^\delta(x,0) = v_0(x) \in H^1.
\endaligned
\ee
With this solution in hand, we need estimates that are uniform in $\delta$ and, clearly, 
the estimates in the proof of Proposition~\ref{prop-400} are not suitable, since they involve 
norms of derivatives of $\varphi^\delta$. Since $\varphi^\delta$ approaches 
a discontinuous function as $\delta \to 0$, these estimates provide no information on the limit $\delta \to 0$.

\begin{lemma}
\label{lem-500}
For $\delta >0$, the solutions $(u^\delta, v^\delta) \in C([0,T]; H^2) \cap C^1( [0,T]; L^2) \times C([0,T]; H^1) \cap C^1([0,T]; L^2)$ 
to the system \eqref{980} satisfy (uniformly in $\delta>0$)
\be
\label{1000}
\aligned
& u^\delta \in L^\infty(0,T; H^1), && u_t^\delta \in L^\infty(0,T; H^{-1}),
\\
& v^\delta \in L^\infty(0,T; H^{-1}), && v_t^\delta \in L^\infty(0,T; H^{-1}),
\\
& u_{xx}^\delta \in L^\infty(0,T; H^{-1}), && v^\delta \varphi^\delta \in L^\infty(0,T; L^2). 
\endaligned
\ee
\end{lemma}

\proof In what follows, the constants $c$ may depend on $t$.
From \eqref{966}, by using the duality $|\langle f, g \rangle | \le \| f \|_{H^{-1}} \|g\|_{H^1}$, we find 
\be
\label{1010}
\aligned
\| u\|_2^2 \le c + c \int_0^T \| u\|_2^2 \,dt + c \, \int_0^T \| v\|_{H^{-1}} \|u\|_{H^1} \,dt
\endaligned
\ee
and, from the energy conservation property \eqref{950}, 
\be
\label{1020}
\aligned
\| u_x\|_2^2 \le c + c \| u\|_2^2  + c \| v\|_{H^{-1}} \|u\|_{H^1} + c \int_{\RR} |\varphi^\delta| v^2 \,dx.
\endaligned
\ee
Now, by multiplying the second equation in \eqref{980} by a test function in $H^1$ and integrating in time, we obtain  
(since $L^2\subset H^{-1}$)
\be
\label{1030}
\aligned
\| v\|^2_{H^{-1}} &\le c + c \, \int_0^T  \int_{\RR} (\varphi^\delta)^2 v^2 \,dx \,dt + c \int_0^T \|u_x\|_2^2 \,dt
\\
&\le c + c \, \int_0^T  \int_{\RR} |\varphi^\delta| v^2 \,dx \,dt + c \int_0^T \|u_x\|_2^2 \,dt.
\endaligned
\ee
The latter inequality is a consequence of $(\varphi^\delta)^2 \le c |\varphi^\delta|$, which follows directly  
from \eqref{150}--\eqref{160}. Also, the same equation yields immediately
\be
\label{1040}
\aligned
\|\del_t v\|^2_{H^{-1}} \le  c \int_{\RR} |\varphi^\delta| v^2 \,dx  + c  \|u_x\|_2^2 . 
\endaligned
\ee

In view of the estimates \eqref{1010}--\eqref{1040}, it is necessary to bound $\int_{\RR} |\varphi^\delta| v^2 \,dx$, only. 
Namely, by multiplying the second equation in \eqref{980} by $v |\varphi^\delta|$, we obtain
\[
\aligned
\frac12\frac d{dt} \int_{\RR} v^2 |\varphi^\delta| \,dx + 2\int_{\RR} (v \varphi^\delta)_x v |\varphi^\delta| \,dx =
\int_{\RR} 2\eps (\Real ru)_x v|\varphi^\delta| \,dx
\endaligned
\]
and, since $\sgn\varphi^\delta(x) = -\sgn x$ (thanks to our choice of regularization), 
\[
\aligned
\int_{\RR} (v \varphi^\delta)_x v |\varphi^\delta| \,dx  &= \frac12 \int_{\RR} (v \varphi^\delta)_x v |\varphi^\delta| \,dx -
\frac12 \int_{\RR} v \varphi^\delta (v |\varphi^\delta|)_x \,dx 
\\
&= \frac12 \int_{-\infty}^0 (v \varphi^\delta)_x v \varphi^\delta \,dx -
\frac12 \int_0^\infty (v \varphi^\delta)_x v \varphi^\delta \,dx 
\\
& \quad - \frac12 \int_{-\infty}^0 v \varphi^\delta (v \varphi^\delta)_x \,dx +
\frac12 \int_0^\infty v \varphi^\delta (v \varphi^\delta)_x \,dx = 0.
\endaligned
\]
Thus, using again $(\varphi^\delta)^2 \le c | \varphi^\delta|$, we arrive at 
\be
\label{1050}
\aligned
\int_{\RR} v^2 | \varphi^\delta| \,dx \le c + c\int_0^T \| u_x\|_2^2 \,dt + c \int_0^T \int_{\RR} v^2 |\varphi^\delta|.
\endaligned
\ee
The estimates \eqref{1010}--\eqref{1050} and Gronwall's lemma lead to the desired properties in \eqref{1000}. The estimates for 
$u_{xx}$ and $u_t$ in \eqref{1000} follow from the equation satisfied by $u$ and
 the proof of Lemma~\ref{lem-500} is now completed.
\endproof

Returning to the proof of Theorem~\ref{thm-200}, we see that \eqref{1000} implies that there exist $u, v, \xi$ such that
(for a subsequence)
\be
\label{1100}
\aligned
&u^\delta \ws u &&\text{ in } L^\infty(0,T; H^1),
\\
&v^\delta \ws v &&\text{ in } L^\infty(0,T; H^{-1}),
\\
&\varphi^\delta v^\delta \ws \xi &&\text{ in } L^\infty(0,T; L^2).
\endaligned
\ee
The distribution $v$ obtained in this way is now checked to be a solution to the second equation in \eqref{200}. As observed in Section~\ref{sec-300}, this is not a trivial question since $v$ has a Dirac measure supported on the line $\Sigma = \{(0,t), t \in(0,T)\},$ and $\varphi $ is discontinuous on $\Sigma$. Following Section~\ref{sec-300}, it is necessary to use the information provided by the equation to define a suitable notion for the  flux $\varphi v$.

Our first observation is as follows: if $\Omega$ is an open set contained in $I \times (0,T)$ with $I \cap \{0\} = \emptyset$, then
$v|_\Omega$ in fact belongs to $ L^\infty(0,T; L^2(\Omega))$ and, moreover, $v$ satisfies the second equation in \eqref{200}
in $\Omega$. To check this, note that, in $\Omega$, the function $\varphi$ is smooth and bounded away from zero.
Therefore, from the last property in \eqref{1100}, we find 
$v^\delta \ws \xi/\varphi$ in $L^\infty(0,T; L^2(\Omega))$
and, along with the second convergence property in \eqref{1100}, we obtain
 $v|_\Omega \in L^\infty(0,T; L^2(\Omega))$ and
$\xi = v\varphi$ almost everywhere in $\Omega$. 
Since, according to \eqref{1100}, $\xi$ satisfies the equation $v_t + 2\xi_x = 2 \varepsilon (r\Real (u))_x$, 
we see that 
\be
\label{1200}
\aligned
v_t + 2(\varphi v)_x = 2 \varepsilon (r\Real (u))_x
\endaligned
\ee
in the sense of distributions in $\Omega$. In other words, the equation \eqref{1200} is satisfied
outside the singular line $\Sigma$.

Given any point $(t,x)$ with $x\neq 0$, one can see using elementary techniques (integrating along the characteristics)
that the equation \eqref{1200} can be uniquely solved with initial data $v_0 \in H^1(\RR)$, since no characteristic 
ever crosses the line $\Sigma$ and the right hand side is an $L^2$ function. 
This provides us with a function $\widetilde v$ defined everywhere except on $\Sigma$, and such that the one--sided limits
$v(0+,t)$ and $v(0-,t)$ exist for all $t>0$.
This shows that $v$ and $\widetilde v$ coincide outside of the singular line $\Sigma$.

Let us now determine the equation satisfied by $\widetilde v$ in the whole of $(0,T) \times \RR$. Let $\phi$ be a test function 
defined on $\RR \times (0,T)$. By computing $\big \langle \widetilde v_t + 2(\varphi \widetilde v)_x , \phi \big \rangle$ in the sense of
distributions and after some easy calculations, we find 
\be
\label{1300}
\widetilde v_t + 2(\varphi \widetilde v)_x = 2 \varepsilon (r\Real (u))_x + J(t) \delta_\Sigma,
\ee
where
\be
\label{1400}
J(t) = 2 \varphi(0+) \widetilde v(0+,t) -  2 \varphi(0-) \widetilde v(0-,t).
\ee
Suppose that $v$ has the form
\be
\label{1450}
\aligned
v = \widetilde v + \Psi(t) \delta_\Sigma, \qquad \Psi \in L^\infty(0,T), 
\endaligned
\ee
that is, $v - \widetilde v$ has the form above and is not a more general distribution.

According to the discussion in Section~\ref{sec-300} and provided we define the distribution product
\be
\label{1500}
\varphi \delta_\Sigma := 0,
\ee
we find (in the sense of distributions)
\[
\aligned
v_t + 2(\varphi v)_x &=  \widetilde v_t + 2(\varphi \widetilde v)_x  
+ \Psi'(t) \delta_\Sigma + 2\varphi \Psi(t)  \delta_\Sigma
\\
& = 2 \varepsilon (r\Real (u))_x + (J(t) + \Psi'(t) )\delta_\Sigma.  
\endaligned
\]
So, $v$ satisfies the equation \eqref{230} under the condition \eqref{1500} if we set $\Psi(t) = -\int_0^t J(s) \,ds$.
Thus, we have established that 
\be
\label{1550}
v = \widetilde v -\int_0^t J(s) \,ds, 
\ee
where $J$ is given by \eqref{1400}. It only remains to show that \eqref{1450} does hold. As seen, $v$ coincides with $\widetilde v$ outside the line $\Sigma$, and so the delicate point is establishing that there exists $\Psi \in L^\infty(0,T)$ satisfying \eqref{1450}. Namely, it is not clear that $v|_\Sigma$ is a bounded function of $t\in(0,T)$, and not a more general distribution.
(One could conceivably have, for instance, $v|_\Sigma = \delta_{t=t_0}$ for some $t_0$.)  
The following technical lemma is in order.

\begin{lemma}
\label{lem-1600}
Let $w \in L^\infty(0,T; H^{-1}(\RR))$ be given such that $\supp w \subset \Sigma = \{ (0,t), t \in (0, T)\}$. Then, there exists $\wbar \in L^\infty(0,T)$ such that $w = \wbar (t)  \delta_\Sigma$.
\end{lemma}

\proof We need to find  a bounded and measurable function $\wbar$ on $(0,T)$ such that
$$
\langle w, \phi (t,x) \rangle  = \int_0^T \wbar (t) \phi(0,t) \,dt \equiv \langle \wbar \otimes \delta_{x=0}, \phi(t,x) \rangle
$$
for every test function $\phi \in \Dcal( (0,T)\times\RR)$ (the space of all compactly supported smooth functions).
First of all, $w$ is a distribution of order at most 1, as this follows from
\[
\aligned
| \langle w ,\phi \rangle | &\le \int_0^T | \langle w(t) , \phi \rangle_{H^{-1}\times H^1} \,dt 
\\
& \le \int_0^T \| w\|_{ L^\infty(0,T; H^{-1})} \| \phi(\cdot,t)\|_{H^1(\RR)} \,dt \le C \| \phi\|_{C^1((0,T) \times \RR)}
\endaligned
\]
with $C$ depending on the support of $\phi$, only.
Therefore, according to H\"ormander \cite[Theorem~2.3.5]{Hormander}, the distribution $w$ has the form
\be
\label{1700}
\aligned
\langle w, \phi \rangle_{\Dcal'((0,T) \times\RR)}  = \langle w_0 , \phi_x (0,t) \rangle_{\Dcal'(0,T)} + \langle w_1, \phi(0,t) \rangle_{\Dcal'(0,T)},
\endaligned
\ee
where $w_0 \in {\Dcal'(0,T)}$ is of order 0 and $w_1 \in {\Dcal'(0,T)}$ is of order at most 1. We will show that, in fact, $w_0 = 0$ and that $w_1$  is actually the bounded function $\wbar$ of interest.

We deal first with the first term in the right-hand side of \eqref{1700} and fix any 
$\phi(t,x)\in \Dcal((0,T) \times \RR)$ such that $\phi(0,t) = 0$ for all $t\in (0,T)$. In that case, \eqref{1700} becomes
\be
\label{1800}
|\langle w, \phi \rangle_{\Dcal'((0,T) \times\RR)} |  = \Big| \int_0^T \phi_x (0,t) \, d w_0(t) \Big|,
\ee
since a distribution of order $0$ is a measure. On the other hand, from the assumption $w \in L^\infty(0,T; H^{-1})$, we have
\be
\label{1900}
|\langle w, \phi \rangle_{\Dcal'((0,T) \times\RR)} |  \le c \| \phi\|_{ L^1(0,T; H^1)}.
\ee
We now construct a family of test functions for which \eqref{1800} and \eqref{1900} lead to a contradiction, unless $w_0 =0$.

Define a family of functions $\{\theta_n(x) \}_{n\ge0}$ in $H^1(\RR)$ by 
\[
\theta_0 (x) =
\left\{
\aligned
& -1 -x, && x \in (-1,-1/2),
\\
& x, && x \in (-1/2, 1/2),
\\
& 1-x, && x \in (1/2, 1),
\\
&0 && \text{otherwise},
\endaligned
\right.
\]
and set $\theta_n(x) := \frac1n \theta_0(nx)$. For all $n$, we have
\be
\label{2000}
\aligned
\theta_n(0) = 0, \qquad \theta_n' (0) = 1, \quad \quad \| \theta_n\|_{H^1(\RR)} \le \frac cn.
\endaligned
\ee
Now let $\phi_n (t,x) = \theta_n (x) \psi(t)$, where $\psi$ is a test function on $(0,T)$. From \eqref{1800} and \eqref{2000}, we get
\be
\label{2100}
\aligned
|\langle w, \phi_n \rangle |  = \Big| \int_0^T \theta_n'(0) \psi(t) \, d w_0(t) \Big|  = \Big| \int_0^T \psi(t) \, d w_0(t) \Big|.
\endaligned
\ee
On the other hand, from \eqref{1900}--\eqref{2000} we find
\be
\label{2200}
\aligned
|\langle w, \phi_n \rangle |  \le \int_0^T | \psi(t) | \frac cn \,dt.
\endaligned
\ee
The estimates \eqref{2100} and \eqref{2200} give 
$\int_0^T \psi(t) dw_0(t) = 0$. 
Since $\psi$ was arbitrary, we conclude that $w_0 =0$.

We now claim that the distribution $w_1$ defined on $(0,T)$ in \eqref{1700} is, in fact, a bounded function. Let $\psi$ be a test function
on $(0,T)$ and let $\phi=\phi(x)$ be a test function on $\RR$ such that $\phi(0) = 1$. Since $w_0 =0$, by using \eqref{1700} we find 
\[
\aligned
|\langle w_1 , \psi \rangle | &= | \langle w , \phi(x) \psi(t) \rangle | 
= \Big| \int_0^T \langle w(t) , \phi(x) \rangle_{H^{-1} \times H^1} \psi(t) \,dt \Big|
\\
& \le  \| w\|_{ L^\infty(0,T; H^{-1})} \| \phi\|_{H^1} \| \psi\|_{L^1(0,T)} = C\| \psi\|_{L^1(0,T)}.
\endaligned
\]
By a density argument, we may suppose that the previous estimate is valid for any $\psi \in L^1(0,T)$. Therefore, $w_1$ is a
continuous functional on $L^1(0,T)$ and so is in $L^\infty(0,T)$. The conclusion of the lemma is precisely \eqref{1700} with $\wbar = w_1$.
\endproof

To complete the  proof of Theorem~\ref{thm-200} we need to establish uniqueness.
Consider two solutions $(u_1,v_1)$ and $(u_2,v_2)$ satisfying the assumptions in Definition~\ref{def-200}. Denote $u$ and $v$ the differences $u_1 - u_2$ and $v_1 - v_2$. Since the system \eqref{200} is linear, $u$ and $v$ satisfy \eqref{200} and so 
uniqueness of solution will follow if we show that any solution with vanishing initial data does vanish on $(0,T) \times \RR$.
For this, the crucial estimate is 
\be
\label{2500}
\aligned
\frac12 \del_t \int_{\RR} | \varphi | \widetilde v^2 \,dx \le 2 \eps \int_{\RR} | \Real (ru)_x | |\varphi \widetilde v|.
\endaligned
\ee
From now, we omit the subscript $\eps$. Suppose that $\widetilde v(0) =0$ and note that $v(0) = \widetilde v(0)$. 
With the estimate \eqref{2500} in hand and from Gronwall's lemma and the fact that $u \in L^\infty(0,T; H^1)$, we find
$\widetilde v(t,x) \equiv 0$. In that case, \eqref{1400} and 
\eqref{1550} imply that, also, $v(t,x) \equiv 0$. 
Then, the first equation becomes $i u_t + u_{xx} = (\varphi + b -2 \eps r^2)u - \eps r^2 \ubar$, 
with $u \in L^\infty(0,T; H^1)$. This allows us to multiply it by $\ubar$ and obtain the conservation of mass \eqref{966} without the latter term. Gronwall's lemma and $u(0) =0$ immediately give us $u\equiv 0$. 

It only remains to check the estimate \eqref{2500}. From the equation \eqref{1300} and after multiplication by $\varphi \widetilde v$ 
and integration over $(-\infty,0)$, we obtain
\[
\aligned
\frac12 \del_t \int_{-\infty}^0 \widetilde v^2 \varphi \,dx + \varphi^2(0-) \widetilde v^2(0-, t) = \int_{-\infty}^0 2 \eps\Real (ru)_x \varphi \widetilde v \,dx
\endaligned
\]
and, since $\varphi \ge 0$ on $(-\infty,0)$,
\[
\aligned
\frac12 \del_t \int_{-\infty}^0 \widetilde v^2 \varphi \,dx \le \int_{-\infty}^0 2 \eps |\Real (ru)_x \varphi \widetilde v| \,dx.
\endaligned
\]
Similarly, integrating over $(0,+\infty)$ and using $\varphi \le 0$, we find 
\[
\aligned
\frac12 \del_t \int^{+\infty}_0 \widetilde v^2 \varphi \,dx \le \int^{+\infty}_0 2 \eps |\Real (ru)_x \varphi \widetilde v| \,dx.
\endaligned
\]
Summing the two previous estimates gives \eqref{2500}. 
(Note in passing that \eqref{2500}, from which uniqueness follows, is satisfied since 
$\varphi$ as an entropy--satisfying shock and would not hold for a ``rarefaction--shock''; see \cite{LeFloch-book} for stability statements about this issue for general solutions to linear hyperbolic equations.)  
This establishes the uniqueness property and completes the proof of Theorem~\ref{thm-200}. \qed


\section{Numerical experiments}
\label{sec-600}

\subsection{Proposed setup}

In this section, we present some numerical computations which provide an illustratation of our theoretical results. We 
compute numerically the solutions of, both, the full nonlinear problem \eqref{100} and the linearized problem \eqref{200}. 
With these simulations in hand, we can check numerically the stability property \eqref{103}. Recall that the linearization approach is based on the formal substitution
$v \mapsto \varphi + \delta v_\delta$ and 
$u \mapsto e^{ibt}r + \delta e^{ibt} u_\delta$, 
where $(u,v)$ are the solutions to the full problem \eqref{100}. This leads us to a linear system satisfied by the perturbation $(u_\delta,v_\delta)$ with initial data $(\ubar_\delta, \vbar_\delta)$.
The claim is that the solution $(u,v)$ to the full system \eqref{100} with perturbed 
initial data $(u_{\mathrm{ref}}(0) + \delta \ubar_\delta, v_{\mathrm{ref}}(0) + \delta \vbar_\delta)$ remains close to the 
perturbed exact solution $(u_{\mathrm{ref}} + \delta e^{ibt} u_\delta, v_{\mathrm{ref}} + \delta v_\delta)$.

Consequently, in order to numerically demonstrate the validity of the linearization procedure, we proceed as follows: 
\begin{enumerate}

\item Compute the solution $(u_\delta,v_\delta)$ to the linearized system \eqref{200} with some initial data $(\ubar_\delta, \vbar_\delta)$.

\item Compute the solution $(u,v)$ to the full system \eqref{100} with perturbed initial data $(u_{\mathrm{ref}}(0) + \delta \ubar_\delta, v_{\mathrm{ref}}(0) + \delta \vbar_\delta)$.

\item Compare $(u,v)$ with the perturbed exact solution $(u_{\mathrm{ref}} + \delta e^{ibt} u_\delta, v_{\mathrm{ref}} + \delta v_\delta)$, and check that the difference remains bounded for any finite time.
\end{enumerate}
\noindent We rely on a finite volume method and will not seek for an in-depth numerical analysis of the solutions, but rather visualize the linearized stability property of interest around our particular reference solution.  


\subsection{Numerical results}

Specifically, the numerical approximations to \eqref{100} are based on the finite volume scheme described in \cite{AmorimFigueira2}. 
The Sch\-r\"o\-din\-ger equation is solved using a finite difference semi-implicit Cranck--Nickolson scheme, and
a Newton algorithm allows us to deal with the
nonlinear term $|u|^2 u$. Burgers equation in \eqref{100} is solved by a semi-implicit Lax--Friedrichs scheme.
Concerning the linearized system \eqref{200}, we employ the Lax--Friedrichs scheme for the second equation and a standard Crank--Nickolson scheme for the Sch\-r\"o\-din\-ger equation. 

Another issue is the numerical computation of the reference solution $r_\eps$ in \eqref{150}. We have employed an explicit Euler scheme after transforming the second--order equation \eqref{150} into a system of two first--order equations. Since the Euler method is not particularly accurate, we rely on a fine mesh (with $20000$ points) so that the solution to \eqref{150} provides a suitable approximation for our purpose. This is confirmed by testing the case $\eps=0$, in which the solution \eqref{175},\eqref{180} is available explicitly.
In our numerical experiments, we choose the parameter values $\eps = 0.1, b=-1.5$ in \eqref{100}, \eqref{150}, and \eqref{160}, and $\delta = 0.1$ in \eqref{102}. The tests are performed with 20000 spatial points on the interval $(-22,22)$ (giving a spatial step $h\simeq 0.002$) and time increment $\tau = 0.0001$. For the perturbations $(\ubar_\delta, \vbar_\delta)$ we choose a gaussian curve $e^{-x^2}$.

\begin{figure}
\includegraphics[width=\linewidth,keepaspectratio=false]{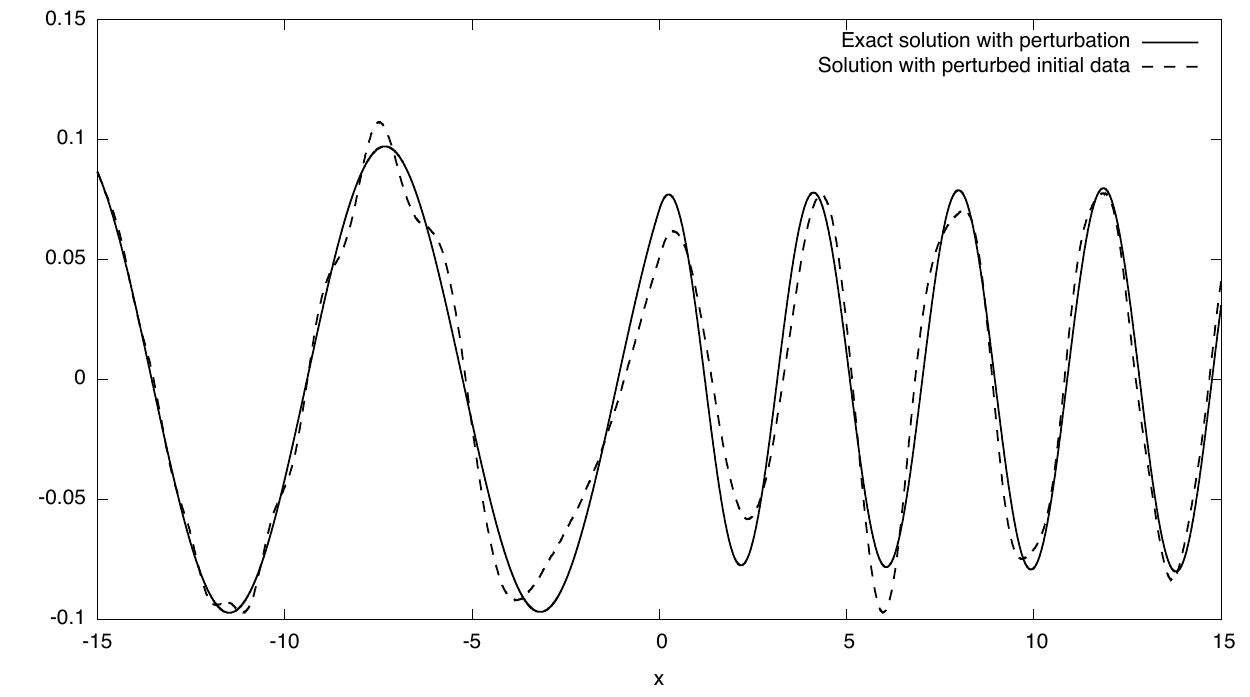}
\caption{Perturbed exact solution and solution with perturbed initial data -- Sch\-r\"o\-din\-ger equation}
\label{fig-10}
\end{figure}

\begin{figure}
\includegraphics[width=\linewidth,keepaspectratio=false]{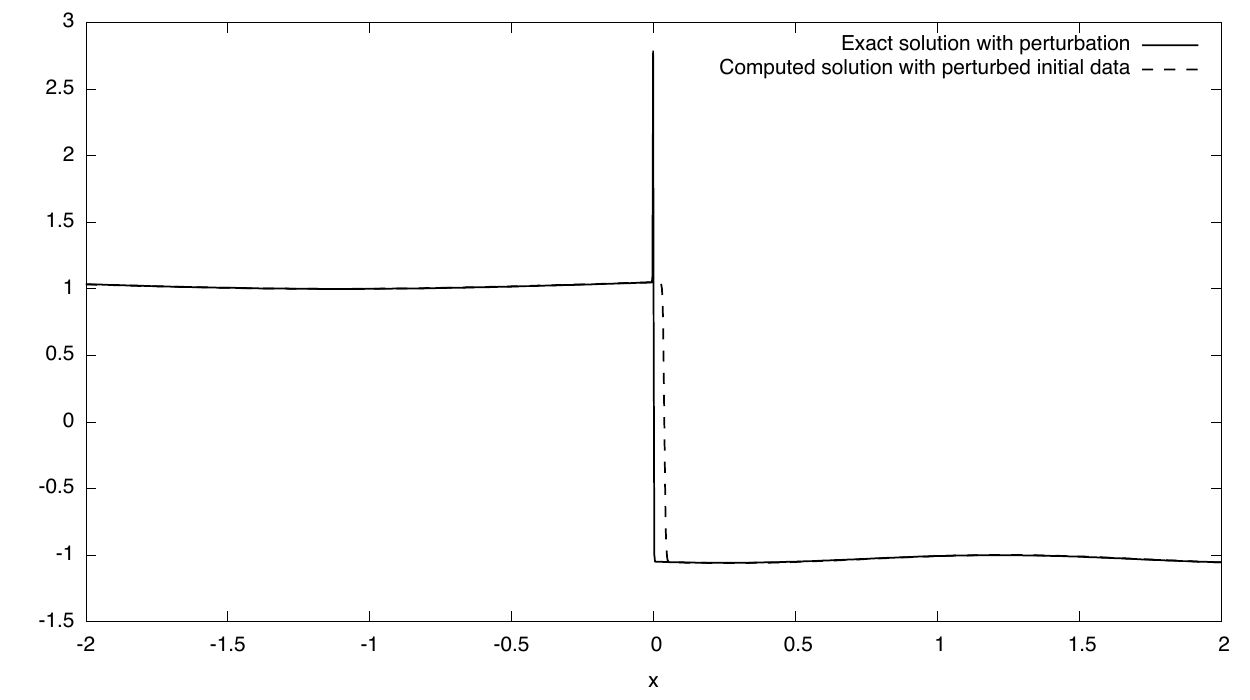}
\caption{Perturbed exact solution and solution with perturbed initial data -- Burgers equation}
\label{fig-20}
\end{figure}

In Figures~\ref{fig-10} and \ref{fig-20}, we present the solution $(u,v)$ to the full system \eqref{100} with perturbed initial data 
$(u_{\mathrm{ref}}(0) + \delta \ubar_\delta, v_{\mathrm{ref}}(0) + \delta \vbar_\delta)$ (dashed line), along with the perturbed exact 
solution $(u_{\mathrm{ref}} + \delta e^{ibt}  u_\delta, v_{\mathrm{ref}} + \delta v_\delta)$ (full line) at the final time $T=1$. 
Recall that the difference between these two quantities should not vanish, but only remain {\sl bounded,} and it is indeed what we observe.

\begin{figure}
\includegraphics[width=\linewidth,keepaspectratio=false]{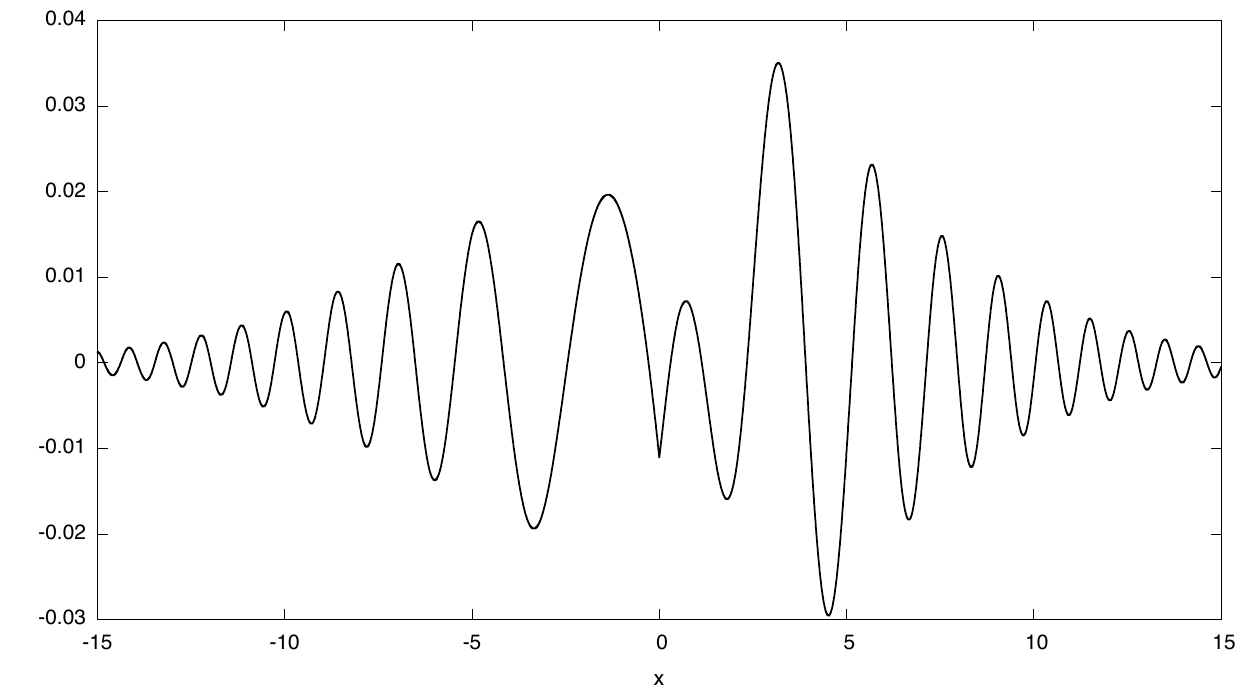}
\caption{Solution to the linearized problem -- Sch\-r\"o\-din\-ger equation}
\label{fig-30}
\end{figure}

\begin{figure}
\includegraphics[width=\linewidth,keepaspectratio=false]{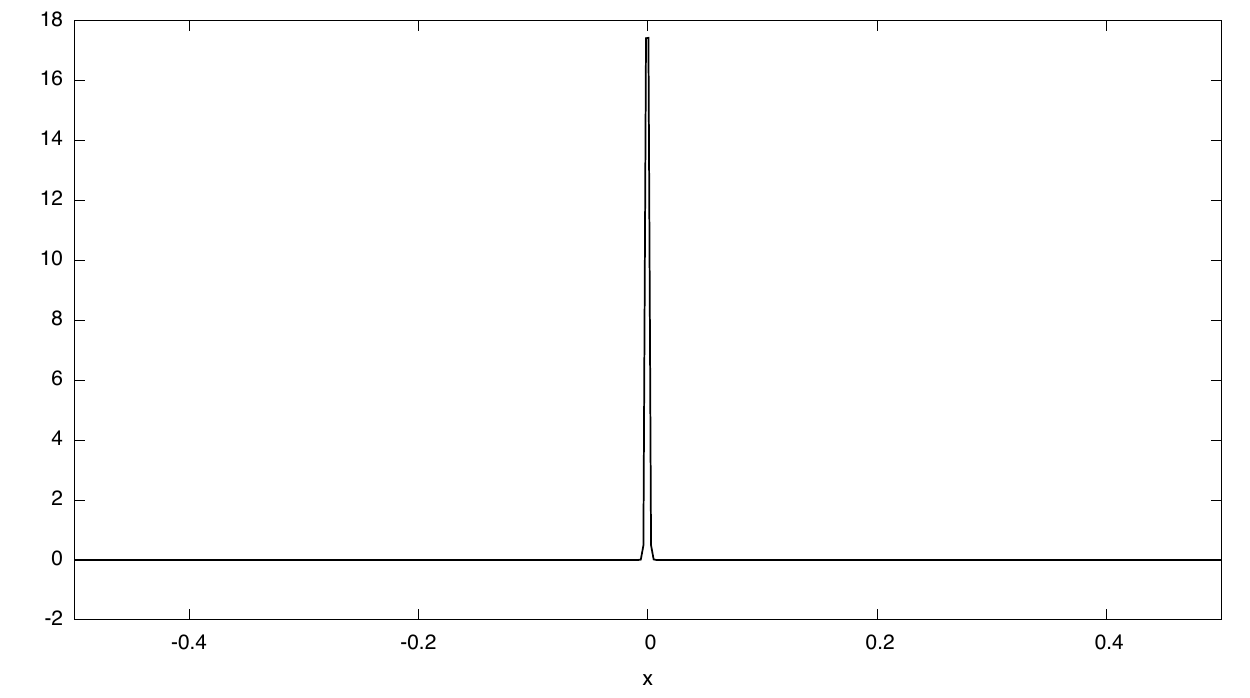}
\caption{Solution to the linearized problem -- Burgers equation.}
\label{fig-40}
\end{figure}

In Figures \ref{fig-30} and \ref{fig-40}, we display the solutions to the linearized system \eqref{200}. 
Recall that the solution to the linear equation with discontinuous coefficient \eqref{200}
contains a Dirac measure on the line $x=0$. Since our scheme contains some numerical viscosity, this measure is smoothed out and appears 
as a steep spike in the solution, which can also be observed in Figure~\ref{fig-20}.

\small 


\section*{Acknowledgements}

The authors are grateful to Luis Sanchez for many discussions. The first three authors were partially supported by the Portuguese Foundation for Science and Technology (FCT) through the grant PTDC/MAT/110613/2009 
and by PEst OE/MAT/UI0209/2011. 
The first author (P.A.) was also supported by the FCT through a \emph{Ci\^encia~2008} fellowship.
The fourth author (PLF) was supported by the Centre National de la Recherche
Scientifique (CNRS) and the Agence Nationale de la Recherche through the grants ANR 2006-2--134423 and ANR SIMI-1-003-01. 


\end{document}